# Consistency and application of moving block bootstrap for non-stationary time series with periodic and almost periodic structure

RAFAL SYNOWIECKI

*Faculty of Applied Mathematics, AGH University of Science and Technology, al. Mickiewicza 30, 30-059 Krakow, Poland. E-mail: rsynowie@agh.edu.pl*

The aim of this paper it to establish sufficient conditions for consistency of moving block bootstrap for non-stationary time series with periodic and almost periodic structure. The parameter of the study is the mean value of the expectation function. Consistency holds in quite general situations: if all joint distributions of the series are periodic, then it suffices to assume the central limit theorem and strong mixing property, together with summability of the autocovariance function. In the case where the mean function is almost periodic, we additionally need uniform boundedness of the fourth moments of the root statistics. It is shown that these theoretical results can be applied in statistical inference concerning the Fourier coefficients of periodically (PC) and almost periodically (APC) correlated time series. A simulation example shows how to use a graphical diagnostic test for significant frequencies and stationarity within these classes of time series.

## 1. Introduction

Moving block bootstrap (MBB), introduced by Künsch [18] and Liu and Singh [22], is a nonparametric bootstrap procedure that can be applied to dependent observations, that is, time series. It consists of calculating the estimator over replicated series obtained by drawing with replacement from the blocks of consecutive data. By means of this resampling procedure, we obtain approximations of unknown distributions of root statistics. Then, on the basis of these approximations, we may construct confidence intervals and test different practical problems. It seems that the MBB procedure is quite well investigated in the case of strictly stationary strong mixing time series (see, e.g., Künsch [18], Radulović [25] and Lahiri [20]) and sufficient conditions for consistency were formulated under quite general conditions for this set-up. As for the non-stationary case, some results have recently been obtained. Fitzenberger [9] and Politis *et al.* [24] give conditions for MBB consistency for the univariate mean without stationarity assumption, whereas the results of Lahiri [19] and Gonçalves and White [13] concern the general heterogeneous time series.

In this paper, we concentrate on a special case of non-stationarity, that is, periodic and almost periodic time series. We develop techniques, presented in Arcones and Giné [2] and Radulović [25], that are based on a general central limit theorem (CLT) for infinitesimal arrays (Araujo and Giné [1]). Our parameter of interest will be the mean value of non-constant, almost periodic mean function $EX_t$, and it is estimated by the estimator $\overline{X}_n = (1/n)\sum_{t=1}^{n} X_t$. Having results concerning this parameter, we can derive results regarding higher-order moments or corresponding Fourier coefficients (see Section 4). Let us argue that the results concerning the univariate mean do not apply here; neither do Theorems 2.2 and 2.1 from Gonçalves and White [13], in which it is required that the time series satisfy the following condition:

$$\frac{1}{n}\sum_{t=1}^{n}\left(EX_t - \frac{1}{n}\sum_{t=1}^{n}EX_t\right)^2 \to 0 \quad \text{for } n \to \infty.$$

However,

$$\frac{1}{n}\sum_{t=1}^{n}\left(EX_t - \frac{1}{n}\sum_{t=1}^{n}EX_t\right)^2 = \frac{1}{n}\sum_{t=1}^{n}(EX_t)^2 - \left(\frac{1}{n}\sum_{t=1}^{n}EX_t\right)^2,$$

and this does not tend to zero if the expectation is a non-constant, almost periodic function (see Lemma A.1 in the Appendix). Finally, the results of Lahiri [19] are very restrictive.

Chan *et al.* [5] and Politis [23] found modifications of the MBB procedure that are applicable to some specific time series with periodic structure. Unfortunately, their procedures require that we know the exact length of the period. Therefore, we cannot use them in the problem of determining significant frequencies. Moreover, these procedures do not generalize to the almost periodic case. Moving block bootstrap itself does not have these restrictions: we do not have to know the period and the procedure can also be easily performed in the almost periodic case.

Statistical inference in the aforementioned special case of non-stationarity is of great practical importance. Many time series data in economics, telecommunications or climatology possess such structure. For examples, we refer the reader to Gardner *el al.* [10] and references therein. Let us add that if we have a series with periodic mean, the period can usually be easily guessed from the plot of the series and we may extract seasonal means to obtain zero-mean time series. However, this cannot be done for higher-order periodicity. Moreover, in the almost periodic case, we are usually unable to extract the mean, which is an almost periodic function. The only solution is to estimate the Fourier coefficients and in order to do this, the significant frequencies should be known. As will be shown in Section 4, the MBB procedure is of considerable help in detecting these significant frequencies.

Section 2 of this paper includes formal definitions of the time series classes that are under consideration. Section 3 presents a description of the MBB algorithm and new results regarding its consistency for strictly periodic and then for almost periodic structured time series. The next section contains the application of MBB in the detection of the significant frequencies for second-order periodically (PC) and almost periodically (APC) correlated time series. All proofs are deferred to the Appendix.

## 2. Classes of time series to be considered

To begin, we give definitions of the time series classes that will be studied in this paper. The order parameter $r$ and the length of the period $T$ are assumed to be positive integers.

**Definition 2.1.** *The time series* $\{X_t : t \in \mathbb{Z}\}$ *is called* strictly periodic of order $r$ (SP($r$)) with period $T$ if, for any $t, \tau_1, \tau_2, \ldots, \tau_{r-1} \in \mathbb{Z}$,

$$(X_t, X_{t+\tau_1}, \ldots, X_{t+\tau_{r-1}}) \stackrel{d}{=} (X_{t+T}, X_{t+\tau_1+T}, \ldots, X_{t+\tau_{r-1}+T}).$$

**Definition 2.2.** *The time series* $\{X_t : t \in \mathbb{Z}\}$ *is called* strictly periodic (SP) *with period* $T$ *if it is* SP($r$) *with period* $T$ *for any* $r \in \mathbb{N}$.

A large class of SP models are time series of the form (Synowiecki [27])

$$X_t = F(Z_t, f(t)),$$

where the time series $\{Z_t\}$ is strictly stationary, the function $f(\cdot)$ is periodic and the function $F(\cdot, f(t))$ is measurable for any $t \in \mathbb{Z}$.

**Definition 2.3.** *The time series* $\{X_t : t \in \mathbb{Z}\}$ *is called* weakly periodic of order $r$ (WP($r$)) with period $T$ if $E|X_t|^r < \infty$ and, for any $t, \tau_1, \tau_2, \ldots, \tau_{r-1} \in \mathbb{Z}$,

$$E(X_t X_{t+\tau_1} \cdots X_{t+\tau_{r-1}}) = E(X_{t+T} X_{t+\tau_1+T} \cdots X_{t+\tau_{r-1}+T}).$$

Before we introduce another definition, we will briefly recall the concept of almost periodic functions. A real- or complex-valued function $f$ is called *almost periodic* if for every $\epsilon > 0$, there exists a number $l_\epsilon$ such that for any interval of length greater than $l_\epsilon$, there exists a number $p_\epsilon$ in this interval such that

$$\sup_{t \in \mathbb{Z}} |f(t + p_\epsilon) - f(t)| < \epsilon$$

(see Corduneanu [6]). Almost periodic functions generalize periodic functions. They possess similar properties, such as boundedness and Fourier representation. As an example, consider the following function:

$$f(t) = \cos(\lambda t).$$

If $t \in \mathbb{R}$, this function is periodic with period $T = 2\pi/\lambda$. However, if $t \in \mathbb{Z}$ and $\lambda \neq 2\pi/m$, where $m \in \mathbb{Z}$, the function $f$ is not periodic, but is almost periodic. The following fact characterizes almost periodic sequences (Corduneanu [6]): the sequence $\{a_n\}$ is almost periodic if and only if there exists an almost periodic function $f$ defined on $\mathbb{R}$ such that $f(n) = a_n$ for all $n \in \mathbb{N}$. The space of almost periodic functions is closed with respect to products, sums and uniform limits. Moreover, for any function $f$ belonging to this class, its mean value, that is, the quantity

$$M_t(f(t)) = \lim_{n \to \infty} \frac{1}{n} \sum_{j=s}^{s+n-1} f(j),$$

exists, uniformly with respect to the number $s$. The subscript $t$ in the symbol $M_t$ is included to emphasize the averaging over the variable $t$.

The next important fact regarding almost periodic functions is that the set

$$\Lambda_f = \{\lambda \in [0, 2\pi) : M_t(f(t)e^{-i\lambda t}) \neq 0\}$$

is countable. If, beyond this, the set $\Lambda_f$ is finite, then the function $f$ has some further desired properties. First, there exists a finite constant $C$, independent of $s$ and $n$, such that

$$\left| \frac{1}{n} \sum_{t=s}^{s+n-1} f(t) - M_t(f(t)) \right| < \frac{C}{n}, \tag{1}$$

which is implied by the discrete counterpart of Lemma 2 from Cambanis *et al.* [4]. Second, the Fourier representation becomes equality (Corduneanu [6]), that is,

$$f(t) = \sum_{\lambda \in \Lambda_f} a(\lambda) e^{-i\lambda t}.$$

This identity is very important in statistical inference of the time series defined below. Note that for a purely periodic function with period $T$, we have

$$\Lambda_f \subset \left\{ \frac{2\pi k}{T} : k = 0, \ldots, T-1 \right\},$$

so this set is always finite. In contradistinction, when the set $\Lambda_f$ is not finite, then it must contain cluster points. Such a situation causes trouble in statistical reasoning (see, e.g., Hurd [14] and Dehay and Leśkow [7]).

**Definition 2.4.** *The time series $\{X_t : t \in \mathbb{Z}\}$ is called* weakly almost periodic of order $r$ (WAP($r$)) *if $E|X_t|^r < \infty$ and, for any $t, \tau_1, \tau_2, \ldots, \tau_{r-1} \in \mathbb{Z}$, the function*

$$E(X_t X_{t+\tau_1} \cdots X_{t+\tau_{r-1}})$$

*is almost periodic in the variable $t$.*

As an example, we could consider amplitude-modulated series of the form $X_t = f(t)Z_t$, which are WAP($r$) provided that the series $\{Z_t : t \in \mathbb{Z}\}$ is WP($r$) and the function $f$ is almost periodic.

It is easy to see that WP($r$) $\subset$ WAP($r$) and, for any two positive integers such that $r_1 < r_2$, we have that SP($r_1$) $\subset$ SP($r_2$) $\subset$ SP. Moreover, each weak periodicity is implied by the corresponding strict periodicity, provided that the appropriate moments exist. Let us add that the class of series that are both WP(1) and WP(2) is identical to the class of periodically correlated (PC) time series in the sense of Gladyshev [12]. A series that is both WAP(1) and WAP(2) is almost periodically correlated (APC) in the sense of Hurd [14]. For statistical inference within these classes, we refer the reader to Hurd and Leśkow [15] and Dehay and Leśkow [7].

## 3. Results regarding MBB for time series with periodic and almost periodic structure

Our aim is to investigate the mean value $\mu = M_t(EX_t)$ of the real-valued time series $\{X_t : t \in \mathbb{Z}\}$ that is (almost) periodic in the first or higher order. The standard estimator obtained from the sample $(X_1, \ldots, X_n)$ is of the form $\overline{X}_n = (1/n) \sum_{t=1}^n X_t$. We wish to determine sufficient conditions for the moving block bootstrap which would enable us to calculate approximations of the quantiles without using the form of asymptotic variance. In order to present the MBB procedure, let us denote the block of $b = b(n)$ consecutive observations as $B_{t,b} = (X_t, \ldots, X_{t+b-1})$ and $k = k(n) = n/b(n)$ which, without loss of generality, is assumed to be integer-valued throughout the whole paper. Let $i_1, i_2, \ldots, i_k$ be i.i.d. random variables with uniform distribution on the set $\{1, 2, \ldots, n-b+1\}$. By joining the blocks $B_{i_1,b}, \ldots, B_{i_k,b}$, we obtain the MBB sample $(X_1^*, \ldots, X_n^*)$ and the MBB version of the estimator takes the form $\overline{X}_n^* = (1/n) \sum_{t=1}^n X_t^*$. In the following, the $P^*$, $E^*$ and Var$^*$ denote quantities obtained from replicated series conditioned on the sample $(X_1, \ldots, X_n)$. The asymptotic variance $\sigma^2$ is always assumed to be positive.

**Theorem 3.1.** *Let $\{X_t : t \in \mathbb{Z}\}$ be a strictly periodic with period $T$, $\alpha$-mixing time series and let $X_t^*$ be generated by the MBB procedure with $b = o(n)$ but $b \to \infty$. Assume that:*

(i) *the autocovariance function is summable, that is,*

$$\sum_{\tau=0}^{\infty} |\mathrm{Cov}(X_t, X_{t+\tau})| < \infty$$

*for all $t = 1, \ldots, T$;*

(ii) *the CLT holds, that is,*

$$\sqrt{n}(\overline{X}_n - \mu) \xrightarrow{d} \mathcal{N}(0, \sigma^2), \qquad (2)$$

*where $\mu = M_t(EX_t)$.*

*Then, MBB is consistent, that is,*

$$\sup_{x \in \mathbb{R}} |P(\sqrt{n}(\overline{X}_n - \mu) < x) - P^*(\sqrt{n}(\overline{X}_n^* - E^*\overline{X}_n^*) < x)| \xrightarrow{P} 0. \tag{3}$$

We can say that the above theorem is a generalization of Theorem 2 of Radulović [25] from the case of strictly stationary to non-stationary strictly periodic time series. We assume periodic structure of all joint distributions and no rate of convergence of $\alpha$-mixing function.

**Corollary 3.1.** *Let $\{X_t : t \in \mathbb{Z}\}$ be a strictly periodic with period $T$, $\alpha$-mixing time series and let $X_t^*$ be generated by the MBB procedure with $b = o(n)$, but $b \to \infty$. Assume that for some $\delta > 0$:*

(i)  $E|X_t|^{2+\delta} < \infty$ for $t = 1, \ldots, T$;
(ii) $\sum_{\tau=1}^{\infty} \alpha_X^{\delta/(2+\delta)}(\tau) < \infty$.

*Then, CLT (2) holds and the MBB procedure is consistent, in the sense of (3).*

The next step will be to formulate a general theorem regarding consistency of moving block bootstrap for non-stationary time series with almost periodic mean.

**Theorem 3.2.** *Let $\{X_t : t \in \mathbb{Z}\}$ be an APC, $\alpha$-mixing time series and let $X_t^*$ be generated by the MBB procedure with $b = o(n)$, but $b \to \infty$. Assume that:*

(i)  *the set $\Lambda = \{\lambda \in [0, 2\pi) : M_t(EX_t e^{-i\lambda t}) \neq 0\}$ is finite;*
(ii) *the autocovariance function is uniformly summable, that is, $|\operatorname{Cov}(X_t, X_{t+\tau})| < c_\tau$, where the sequence $\{c_\tau\}_{\tau=0}^{\infty}$ is summable;*
(iii) *there exists a finite constant $K$ that does not depend on $b = b(n) \leq n$ and $n$, such that*

$$\sup_{s=1,\ldots,n-b+1} E\left(\frac{1}{\sqrt{b}} \sum_{t=s}^{s+b-1} (X_t - EX_t)\right)^4 \leq K; \tag{4}$$

(iv) *the CLT holds, that is,*

$$\sqrt{n}(\overline{X}_n - \mu) \xrightarrow{d} \mathcal{N}(0, \sigma^2), \tag{5}$$

*where $\mu = M_t(EX_t)$.*

*Then, the MBB is consistent, that is,*

$$\sup_{x \in \mathbb{R}} |P(\sqrt{n}(\overline{X}_n - \mu) < x) - P^*(\sqrt{n}(\overline{X}_n^* - E^*\overline{X}_n^*) < x)| \xrightarrow{P} 0. \tag{6}$$

We will now describe two specific situations in which Theorem 3.2 is satisfied. Assumptions (ii) and (iii) can be guaranteed by an appropriate mixing rate and uniform boundedness of moments of the series.

**Corollary 3.2.** *Let $\{X_t : t \in \mathbb{Z}\}$ be an APC, $\alpha$-mixing time series and let $X_t^*$ be generated by the MBB procedure with $b = o(n)$, but $b \to \infty$. Assume that:*

  (i) *the set $\Lambda = \{\lambda \in [0, 2\pi) : M_t(EX_t \mathrm{e}^{-\mathrm{i}\lambda t}) \neq 0\}$ is finite;*
  (ii) *$\sup_{t \in \mathbb{Z}} E|X_t|^{4+\delta} < \infty$ for some $\delta > 0$;*
  (iii) *$\sum_{\tau=1}^{\infty} \tau \alpha_X^{\delta/(4+\delta)}(\tau) < \infty$;*
  (iv) *the CLT holds (i.e., (5) is satisfied).*

*Then, the MBB procedure is consistent, in the sense of (6).*

First, note that (5) can be obtained by, for example, Theorem B.0.1 of Politis *et al.* [24]. We can easily verify its condition (B.2) by means of Lemma A.6. A slightly stronger mixing rate should be assumed, that is,

$$\sum_{\tau=0}^{\infty} (\tau+1)^2 \alpha_X^{(2+\delta)/(10+\delta)}(\tau) < \infty,$$

where the parameter $\delta$ is defined in Corollary 3.2. Second, we can take a constant expectation function and the series for which, instead of WAP(2) assumption,

$$\sup_{s=1,\ldots,n-b+1} \mathrm{Var}\left(\frac{1}{\sqrt{b}} \sum_{t=s}^{s+b-1} X_t\right) \to \sigma^2 \qquad \text{for } n \to \infty.$$

Then, Corollary 3.2 can be viewed as a generalization of Theorem 4.4.2 (the univariate mean case) from Politis *et al.* [24] since the assumption regarding the mixing rate is weaker.

The last result of this section concerns uniformly bounded time series.

**Corollary 3.3.** *Let $\{X_t : t \in \mathbb{Z}\}$ be an APC, $\alpha$-mixing time series and let $X_t^*$ be generated by the MBB procedure with $b = o(n)$, but $b \to \infty$. Assume that:*

  (i) *the set $\Lambda = \{\lambda \in [0, 2\pi) : M_t(EX_t \mathrm{e}^{-\mathrm{i}\lambda t}) \neq 0\}$ is finite;*
  (ii) *the series $\{X_t\}$ is a.s. uniformly bounded;*
  (iii) *$\alpha_X(\tau) = O(\tau^{-2})$;*
  (iii) *the CLT holds (i.e., (5) is satisfied).*

*Then the MBB procedure is consistent, in the sense of (6).*

Therefore, in the case of uniformly bounded random variables, further relaxation regarding the mixing rate can be allowed.

## 4. Application to the Fourier coefficients of the autocovariance function

In the case of periodically (PC) and almost periodically (APC) correlated time series, statistical inference focuses mainly on second-order properties of the series. If no assumption regarding the APC model is made, it is based on the Fourier representation of the function $B(t,\tau) = E(X_t X_{t+\tau})$ (Hurd and Leśkow [15], Dehay and Leśkow [7]), that is

$$B(t,\tau) = \sum_{\lambda \in \Lambda_\tau} a(\lambda,\tau) e^{i\lambda t}.$$

We let

$$a(\lambda,\tau) = M_t(B(t,\tau) e^{-i\lambda t})$$

and the set $\Lambda_\tau = \{\lambda : a(\lambda,\tau) \neq 0\}$ is assumed to be finite. In the engineering literature, $a(\lambda,\tau)$ is called a *cyclic autocorrelation function*, while the elements of the set $\Lambda_\tau$ are called *cyclic frequencies*.

If we observe the real-valued series $\{X_t : t \in \mathbb{Z}\}$ for $t = 1,\ldots,n$, the estimator of the parameter $a(\lambda,\tau)$ is of the form

$$\hat{a}_n(\lambda,\tau) = \frac{1}{n-|\tau|} \sum_{t=1-\min\{\tau,0\}}^{n-\max\{\tau,0\}} X_t X_{t+\tau} e^{-i\lambda t}.$$

Its properties, such as strong consistency and asymptotic normality, were studied for the case of $\varphi$-mixing stochastic processes in Hurd and Leśkow [15] and Dehay and Leśkow [7]. Recently, these results have been generalized to the case of $\alpha$-mixing APC time series (Leśkow and Synowiecki [21]). The aforementioned sufficient conditions for the CLT for the estimator $\hat{a}_n(\lambda,\tau)$ require the mixing rate

$$\sum_{\tau=1}^{\infty} (\tau+1)^2 \alpha_X^{\delta/(4+\delta)}(\tau) < \infty,$$

provided that $\sup_t E|X_t|^{4+4\delta} < \infty$ and that the series is also WAP(4). Unfortunately, the form of asymptotic variance is quite complicated and depends on fourth-order mean values, so it is not applicable in practice. This justifies why it is important to develop resampling methods such as MBB within the class of PC and APC time series.

**Remark 4.1.** Without loss of generality, assume that $\tau \geq 0$. For the APC time series $\{X_t : t \in \mathbb{Z}\}$, define the time series

$$W_t(\lambda,\tau) = X_t X_{t+\tau} e^{-i\lambda t}.$$

Its expectation function is almost periodic and, in fact, for all $\lambda \in [0, 2\pi)$, apart from at most a finite number of them, not purely periodic. It is easy to see that

$$M_t(EW_t(\lambda, \tau)) = \lim_{n \to \infty} \frac{1}{n} \sum_{t=1}^{n} E(X_{t+\tau} X_t e^{-i\lambda t}) = a(\lambda, \tau)$$

and

$$\overline{W}_{n-\tau}(\lambda, \tau) = \frac{1}{n - \tau} \sum_{t=1}^{n-\tau} W_t = \hat{a}_n(\lambda, \tau).$$

Therefore, it will be possible of use the results regarding MBB for the expectation function of the estimator $\hat{a}_n(\lambda, \tau)$. Let $(W_1^*(\lambda, \tau), \ldots, W_{n-\tau}^*(\lambda, \tau))$ be an MBB sample obtained from the sample

$$(W_1(\lambda, \tau), \ldots, W_{n-\tau}(\lambda, \tau)) = (X_1 X_{1+\tau} e^{-i\lambda}, \ldots, X_{n-\tau} X_n e^{-i\lambda(n-\tau)}).$$

The MBB version of the estimator $\hat{a}_n(\lambda, \tau)$ is defined as

$$\hat{a}_n^*(\lambda, \tau) = \overline{W}_{n-\tau}^*(\lambda, \tau) = \frac{1}{n - \tau} \sum_{t=1}^{n-\tau} W_t^*(\lambda, \tau).$$

We now provide a theorem regarding the consistency of the MBB estimator $\hat{a}_n^*(\lambda, \tau)$. The inequalities between complex numbers involving $a(\lambda, \tau)$, $\hat{a}_n(\lambda, \tau)$ and $(x, y)$ are to be understood componentwise.

**Theorem 4.1.** *Let $\{X_t : t \in \mathbb{Z}\}$ be an APC and WAP(4) time series that satisfies the following conditions:*

  (i) *the set $\Lambda_\tau = \{\lambda \in [0, 2\pi) : M_t(B(t, \tau)e^{-i\lambda t}) \neq 0\}$ is finite;*
  (ii) $\sup_{t \in \mathbb{Z}} E|X_t|^{8+2\delta} < \infty$;
  (iii) $\sum_{k=1}^{\infty} \tau \alpha_X^{\delta/(4+\delta)}(\tau) < \infty$;
  (iv) *the CLT for the estimator $\hat{a}_n(\lambda, \tau)$ holds, that is,*

$$\sqrt{n}(\hat{a}_n(\lambda, \tau) - a(\lambda, \tau)) \xrightarrow{d} \mathcal{N}_2(0, \Sigma) \qquad \text{for } n \to \infty$$

*and $\det(\Sigma) \neq 0$.*

*Then, MBB with $b = o(n)$, but $b \to \infty$, for the estimator $\hat{a}_n(\lambda, \tau)$ is consistent, that is,*

$$\sup_{(x,y) \in \mathbb{R}^2} |P(\sqrt{n}(\hat{a}_n(\lambda, \tau) - a(\lambda, \tau)) < (x, y))$$
$$- P^*(\sqrt{n}(\hat{a}_n^*(\lambda, \tau) - E^*(\hat{a}_n^*(\lambda, \tau)))) < (x, y))| \xrightarrow{P} 0.$$

The MBB procedure for the estimator $\hat{a}_n(\lambda, \tau)$ that is described above may be used in the problem of determining significant frequencies (i.e., the frequencies that belong to the set $\Lambda$), which, in the PC case, is equivalent to the identification of the period. In practical applications (see, e.g., Yeung and Gardner [28], Dehay and Leśkow [8] and Gardner *et al.* [10]), one calculates the values of the estimator $\hat{a}_n(\lambda, \tau)$ for $\lambda \in [0, 2\pi)$. The frequencies for which spikes are obtained are then chosen to be significant. Unfortunately, this choice is made arbitrarily, as it is not possible to construct reasonable confidence intervals based on the asymptotic distribution; see Dehay and Leśkow [7] and Leśkow and Synowiecki [21] for the exact form of the asymptotic variance. Therefore, we propose to use MBB, which provides an easy way to obtain the pointwise consistent confidence intervals. Certainly, it would be most desirable to conduct simultaneous inference with respect to $\lambda$ for $a(\lambda, \tau)$. The research of Dehay and Leśkow [7] presents a functional CLT in $\tau$ for the root statistics. However, thus far, no such result exists in the argument $\lambda$.

## 5. Simulation example

We will use the series simulated from the PAR(1) model

$$X_t = a_t X_{t-1} + \epsilon_t,$$

where

$$a_t = \frac{2}{3} + \frac{1}{3}\sin\left(\frac{2\pi t}{3}\right)$$

and $\epsilon_1, \epsilon_2, \ldots$ are i.i.d. with the standard normal distribution. This model is PC with period $T = 3$ (see Bloomfield *et al.* [3]) and $\alpha$-mixing with a geometrically decaying mixing function. Figure 1(a) depicts the series itself. Note that it is not possible to guess the second-order period from this plot. Figure 1(b) shows the values of the real part of the estimator $\hat{a}_n(\lambda, 1)$ with respect to $\lambda \in [0, \pi]$. Since $\hat{a}_n(2\pi - \lambda, \tau) = \overline{\hat{a}_n(\lambda, \tau)}$, we do not have to plot the values for the interval $(\pi, 2\pi)$. Finally, Figure 1(c) corresponds to the imaginary part of the estimator. At this stage, we are unable to tell which frequencies are significant. We can observe several spikes, none of which seems to dominate.

Consider the following testing problem:

$$H_0: a(\lambda, \tau) = 0;$$
$$H_1: a(\lambda, \tau) \neq 0.$$

To verify $H_0$, we will use $T_n(\lambda, \tau) = \sqrt{n}\,\hat{a}_n(\lambda, \tau)$, which, under $H_0$, has zero-mean asymptotic normal distribution. The quantiles of $T_n(\lambda, \tau)$ are obtained by means of the MBB procedure, that is, we approximate the quantiles of $T_n^*(\lambda, \tau)$ via Monte Carlo simulation. We used the block size $b = 30$. Figure 2 depicts the values of the estimator $\hat{a}_n(\lambda, 1)$ (dotted lines) compared to this quantiles (rescaled by $\sqrt{n}$) of orders 0.05 and 0.95 (solid lines). The fact that the dotted lines cross the solid line means that, at this point, we reject the hypothesis $H_0$. The only points at which this happens are $\lambda \approx 2.1 \approx 2\pi/3$ and $\lambda \approx 0$. Therefore, we conclude that our series has period $T = 3$ because, in this case, $\Lambda_\tau \subset \{0, 2\pi/3, 4\pi/3\}$. Let us add that if there were not any significant spikes for any $\tau$, we might conclude that the series is stationary. Therefore, this graphical test can also be viewed as a stationarity test within the class of APC time series.

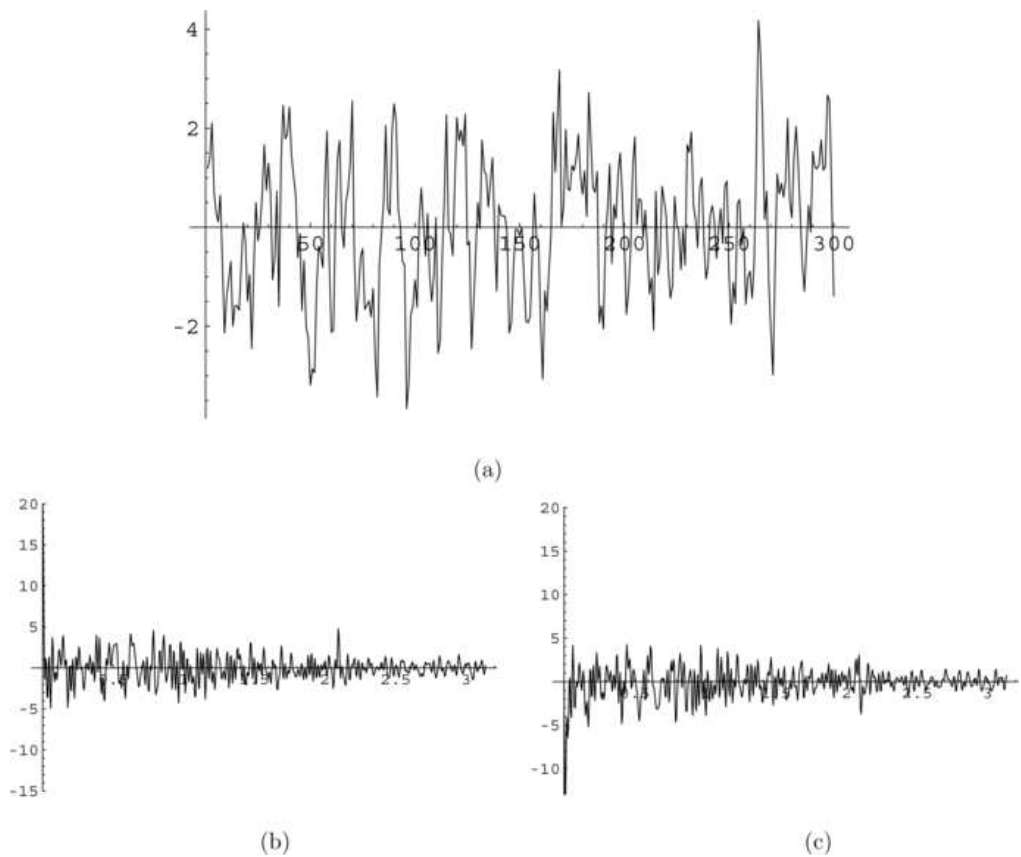

**Figure 1.** Unsuccessful detection of significant frequencies for the simulated series. (a) The PAR(1) series, the sample size here being $n = 300$. (b) The values of $\operatorname{Re}(\hat{a}_n(\lambda, 1))$ with respect to $0 \leq \lambda \leq \pi$. (c) The values of $\operatorname{Im}(\hat{a}_n(\lambda, 1))$ with respect to $0 \leq \lambda \leq \pi$.

## 6. Conclusions

The aim of this paper was to show that the classical MBB procedure works well for those non-stationary time series that have periodic or almost periodic structure. For the time series that have all joint distributions periodic, consistency holds under very general conditions; for the almost periodic case, we need some more stringent assumptions regarding moments of the series and the mixing rate. As for further research, it would be very desirable to estimate within these classes of models the rate of optimal block size $b$, as well as to compare (theoretically) performance of MBB with that of the procedures proposed by Chan *et al.* [5] and Politis [23] for the case of periodic series.

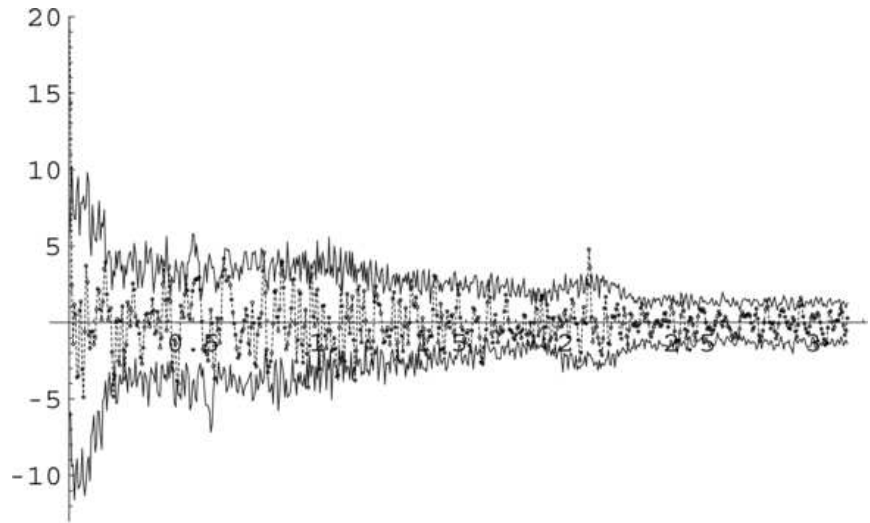

(a)

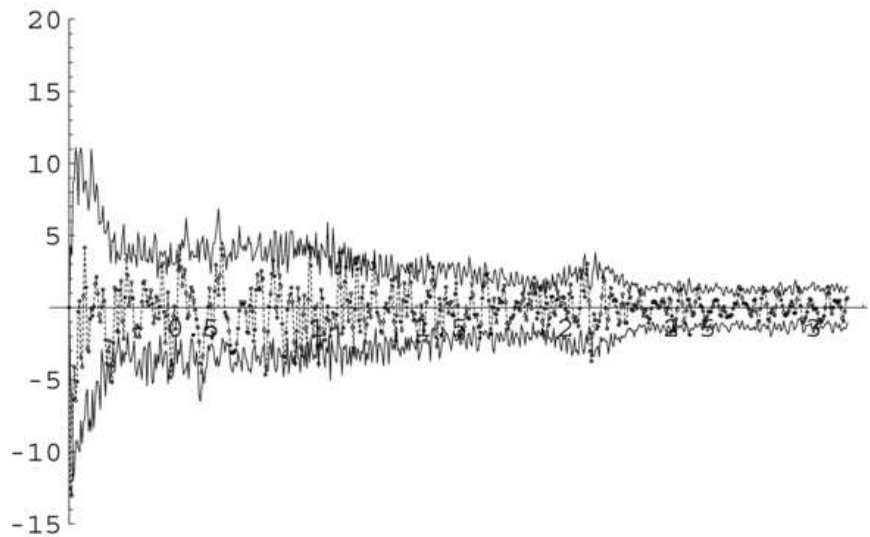

(b)

**Figure 2.** Successful MBB-based detection of significant frequencies for the simulated series. (a) The values of $\text{Re}(\hat{a}_n(\lambda, 1))$ (dotted line) and MBB confidence intervals (solid lines) with respect to $0 \leq \lambda \leq \pi$. (b) The values of $\text{Im}(\hat{a}_n(\lambda, 1))$ (dotted line) and MBB confidence intervals (solid lines) with respect to $0 \leq \lambda \leq \pi$.

# Appendix

**Lemma A.1.** *Let the real-valued function $f$ be almost periodic. Assume that the set $\Lambda_f = \{\lambda \in [0, 2\pi) : M_t(f(t)e^{-i\lambda t}) \neq 0\}$ is finite. If we have*

$$M_t(f^2(t)) = M_t^2(f(t)),$$

*then the function $f$ is constant.*

**Proof.** Using the Fourier representation of almost periodic functions, we have that

$$f(t) = \sum_{\lambda \in \Lambda_f} a(\lambda) e^{i\lambda t}.$$

Therefore,

$$M_t(f^2(t)) = \sum_{\lambda \in \Lambda_f} a(\lambda) a(2\pi - \lambda).$$

Since $a(2\pi - \lambda) = \overline{a(\lambda)}$, the assumption implies that $a^2(0) = \sum_{\lambda \in \Lambda_f} |a(\lambda)|^2$, so $\Lambda_f = \{0\}$. □

**Lemma A.2.** *Let $(X_1, \ldots, X_n)$ be a sample from the time series $\{X_t : t \in \mathbb{Z}\}$ that is strictly periodic with period $T$. Then, the triangular array $\{Y_{n,t} : t = 1, \ldots, n - b + 1\}$, where $Y_{n,t} = X_t + \cdots + X_{t+b-1}$, is row-wise strictly periodic with the same period $T$ for any sequence of positive integers $b = b(n) \leq n$.*

**Proof.** Since, for any $r \in \mathbb{N}$, $t, \tau_1, \ldots, \tau_{r-1} \in \mathbb{Z}$,

$$(X_t, X_{t+1}, \ldots, X_{t+b-1}, X_{t+\tau_1}, X_{t+\tau_1+1}, \ldots, X_{t+\tau_1+b-1},$$
$$\ldots, X_{t+\tau_{r-1}}, X_{t+\tau_{r-1}+1}, \ldots, X_{t+\tau_{r-1}+b-1})$$
$$\stackrel{d}{=} (X_{t+T}, \ldots, X_{t+T+b-1}, X_{t+\tau_1+T}, \ldots, X_{t+\tau_1+T+b-1},$$
$$\ldots, X_{t+\tau_{r-1}+T}, \ldots, X_{t+\tau_{r-1}+T+b-1}),$$

calculating the values of these vectors under the Borel measurable mapping that sums $b$ successive elements of an $rb$-dimensional vector, we obtain

$$(Y_{n,t}, Y_{n,t+\tau_1}, \ldots, Y_{n,t+\tau_{r-1}}) \stackrel{d}{=} (Y_{n,t+T}, Y_{n,t+\tau_1+T}, \ldots, Y_{n,t+\tau_{r-1}+T}). \qquad \square$$

**Lemma A.3.** *If the time series $\{X_t : t \in \mathbb{Z}\}$ is strictly periodic (or $\mathrm{SP}(r)$) with period $T$ and the function $f : \mathbb{R} \longrightarrow \mathbb{R}$ is Borel measurable, then the time series $\{f(X_t) : t \in \mathbb{Z}\}$ is also strictly periodic (or $\mathrm{SP}(r)$) with the same period $T$.*

**Proof.** For any $r \in \mathbb{N}$, define the function $g : \mathbb{R}^r \longrightarrow \mathbb{R}^r$ as

$$g(x_1, \ldots, x_r) = (f(x_1), \ldots, f(x_r)).$$

Knowing that for any $t, \tau_1, \ldots, \tau_{r-1} \in \mathbb{Z}$,

$$(X_t, X_{t+\tau_1}, \ldots, X_{t+\tau_{r-1}}) \stackrel{d}{=} (X_{t+T}, X_{t+\tau_1+T}, \ldots, X_{t+\tau_{r-1}+T}),$$

calculating the values of these vectors under the measurable mapping $g$, we obtain

$$(f(X_t), f(X_{t+\tau_1}), \ldots, f(X_{t+\tau_{r-1}})) \stackrel{d}{=} (f(X_{t+T}), f(X_{t+\tau_1+T}), \ldots, f(X_{t+\tau_{r-1}+T})).$$

$\square$

**Lemma A.4.** *Let $\{X_{n,t} : t = 1, \ldots, d_n\}$, where $d_n \to \infty$, be a triangular array of real random variables, which is row-wise $\mathrm{SP}(1)$ with the same period $T$. Assume that:*

(i) *for $t = 1, \ldots, T$, the series $\{X_{n,t}\}_{n=n_t}^{\infty}$ are uniformly integrable (we denote by $X_{n,n_t}$ the first element in column $t$ of the array $\{X_{n,t}\}$);*

(ii) *the following limits exist and are finite*

$$\mu_1 = \lim_{n \to \infty} EX_{n,1},$$
$$\mu_2 = \lim_{n \to \infty} EX_{n,2},$$
$$\vdots$$
$$\mu_T = \lim_{n \to \infty} EX_{n,T};$$

(iii) *there exists a triangular array of non-negative real numbers $\{a_{n,\tau} : \tau = 0, \ldots, d_n - 1\}$ such that*

$$\frac{1}{d_n} \sum_{\tau=0}^{d_n-1} a_{n,\tau} \to 0 \quad \text{for } n \to \infty$$

*and, for every $A \in \mathbb{R}$, $t = 1, \ldots, T$ and $\tau = 0, 1, \ldots, d_n - t$,*

$$\mathrm{Cov}(X_{n,t}\mathbf{1}_{|X_{n,t}|<A}, X_{n,t+\tau}\mathbf{1}_{|X_{n,t+\tau}|<A}) \leq A^2 a_{n,\tau}.$$

*We then have*

$$\frac{1}{d_n} \sum_{t=1}^{d_n} X_{n,t} \xrightarrow{P} \mu,$$

where $\mu = (1/T)\sum_{t=1}^{T} \mu_t$.

**Proof.** The technique to be used here is a modification of the proof of Lemma 1 from Radulović [25]. Let $Y_{n,t} = X_{n,t} - EX_{n,t}$. Due to the periodicity of $EX_{n,t}$, it suffices to prove that

$$\frac{1}{d_n}\sum_{t=1}^{d_n} Y_{n,t} \xrightarrow{P} 0.$$

Letting $\widetilde{Y}_{n,t} = Y_{n,t}\mathbf{1}_{|Y_{n,t}|<A_n}$ and $\overline{Y}_{n,t} = Y_{n,t}\mathbf{1}_{|Y_{n,t}|\geq A_n}$, the above convergence follows from

$$\frac{1}{d_n}\sum_{t=1}^{d_n}(\widetilde{Y}_{n,t} - E\widetilde{Y}_{n,t}) \xrightarrow{P} 0 \tag{7}$$

and

$$\frac{1}{d_n}\sum_{t=1}^{d_n}(\overline{Y}_{n,t} - E\overline{Y}_{n,t}) \xrightarrow{P} 0. \tag{8}$$

To prove (7), note that, for any $\epsilon > 0$,

$$P\left(\left|\frac{1}{d_n}\sum_{t=1}^{d_n}(\widetilde{Y}_{n,t} - E\widetilde{Y}_{n,t})\right| > \epsilon\right) \leq \frac{1}{d_n^2\epsilon^2}\operatorname{Var}\left(\sum_{t=1}^{d_n}(\widetilde{Y}_{n,t} - E\widetilde{Y}_{n,t})\right)$$

$$\leq \frac{2}{d_n^2\epsilon^2}\sum_{t=1}^{d_n}\sum_{\tau=0}^{d_n-t}|\operatorname{Cov}(\widetilde{Y}_{n,t},\widetilde{Y}_{n,t+\tau})|$$

$$\leq \frac{2}{d_n^2\epsilon^2}\sum_{t=1}^{d_n}\sum_{\tau=0}^{d_n-1}A_n^2 a_{n,\tau} = \frac{2}{d_n\epsilon^2}\sum_{\tau=0}^{d_n-1}A_n^2 a_{n,\tau}.$$

Letting

$$A_n = \left(\frac{1}{d_n}\sum_{\tau=0}^{d_n-1} a_{n,\tau}\right)^{-1/4},$$

it is easy to see that $A_n \to \infty$ and we therefore obtain that

$$P\left(\left|\frac{1}{d_n}\sum_{t=1}^{d_n}(\widetilde{Y}_{n,t} - E\widetilde{Y}_{n,t})\right| > \epsilon\right) \leq \frac{2}{\epsilon^2}\left(\frac{1}{d_n}\sum_{\tau=0}^{d_n-1} a_{n,\tau}\right)^{1/2} \to 0 \quad \text{for } n \to \infty.$$

In order to prove (8), note that, by Lemma A.3 with the function $f(x) = x\mathbf{1}_{|x|>A}$, we have, for any $\epsilon > 0$,

$$P\left(\left|\frac{1}{d_n}\sum_{t=1}^{d_n}(\overline{Y}_{n,t} - E\overline{Y}_{n,t})\right| > \epsilon\right) \leq \frac{1}{d_n\epsilon}E\left|\sum_{t=1}^{d_n}(\overline{Y}_{n,t} - E\overline{Y}_{n,t})\right|$$

$$\leq \frac{2}{d_n\epsilon}\sum_{t=1}^{d_n}E|\overline{Y}_{n,t}|$$

$$\leq \frac{2(\lfloor d_n/T \rfloor + 1)}{d_n\epsilon}\sum_{t=1}^{T}E|X_{n,t}|\mathbf{1}_{|X_{n,t}|\geq A_n}.$$

By uniform integrability of the sequences $\{X_{n,t}\}_{n=n_t}^{\infty}$ for $t = 1,\ldots,T$, we obtain the desired convergence to zero. $\square$

**Lemma A.5.** *Let $(X_1,\ldots,X_n)$ be a sample from the time series $\{X_t : t \in \mathbb{Z}\}$ that is WAP(1). Assume that*

(i) *the set $\Lambda = \{\lambda : M_t(EX_t e^{-i\lambda t}) \neq 0\}$ is finite;*

(ii) *there exists a finite constant $K$ that does not depend on $b = b(n) \leq n$ and $n$, such that*

$$\sup_{s=1,\ldots,n-b+1} E\left(\frac{1}{\sqrt{b}}\sum_{t=s}^{s+b-1}(X_t - EX_t)\right)^4 < K.$$

*There then exists a finite constant $K'$ that does not depend on $b$ and $n$, such that*

$$\sup_{s=1,\ldots,n-b+1} E\left(\frac{1}{\sqrt{b}}\sum_{t=s}^{s+b-1}(X_t - \mu)\right)^4 < K',$$

*where $\mu = M_t(EX_t)$.*

**Proof.** We have

$$E\left(\frac{1}{\sqrt{b}}\sum_{t=s}^{s+b-1}(X_t - \mu)\right)^4 = E\left(\frac{1}{\sqrt{b}}\sum_{t=s}^{s+b-1}(X_t - EX_t) + \frac{1}{\sqrt{b}}\sum_{t=s}^{s+b-1}(EX_t - \mu)\right)^4$$

$$= E\left(\frac{1}{\sqrt{b}}\sum_{t=s}^{s+b-1}(X_t - EX_t)\right)^4$$

$$+ 4\left(\frac{1}{\sqrt{b}}\sum_{t=s}^{s+b-1}(EX_t - \mu)\right)E\left(\frac{1}{\sqrt{b}}\sum_{t=s}^{s+b-1}(X_t - EX_t)\right)^3$$

$$+ 6\left(\frac{1}{\sqrt{b}} \sum_{t=s}^{s+b-1} (EX_t - \mu)\right)^2 E\left(\frac{1}{\sqrt{b}} \sum_{t=s}^{s+b-1} (X_t - EX_t)\right)^2$$

$$+ 4\left(\frac{1}{\sqrt{b}} \sum_{t=s}^{s+b-1} (EX_t - \mu)\right)^3 E\left(\frac{1}{\sqrt{b}} \sum_{t=s}^{s+b-1} (X_t - EX_t)\right)$$

$$+ \left(\frac{1}{\sqrt{b}} \sum_{t=s}^{s+b-1} (EX_t - \mu)\right)^4.$$

Due to (1), we have

$$\frac{1}{\sqrt{b}} \sum_{t=s}^{s+b-1} (EX_t - \mu) = O\left(\frac{1}{\sqrt{b}}\right).$$

Therefore,

$$E\left(\frac{1}{\sqrt{b}} \sum_{t=s}^{s+b-1} (X_t - \mu)\right)^4$$

$$= E\left(\frac{1}{\sqrt{b}} \sum_{t=s}^{s+b-1} (X_t - EX_t)\right)^4 + O\left(\frac{1}{\sqrt{b}}\right) E\left(\frac{1}{\sqrt{b}} \sum_{t=s}^{s+b-1} (X_t - EX_t)\right)^3$$

$$+ O\left(\frac{1}{b}\right) E\left(\frac{1}{\sqrt{b}} \sum_{t=s}^{s+b-1} (X_t - EX_t)\right)^2 + O\left(\frac{1}{b^2}\right).$$

The first term is uniformly bounded, by assumption (ii). For the second and third terms, we have

$$\left|E\left(\frac{1}{\sqrt{b}} \sum_{t=s}^{s+b-1} (X_t - EX_t)\right)^2\right| \leq \left\{E\left(\frac{1}{\sqrt{b}} \sum_{t=s}^{s+b-1} (X_t - EX_t)\right)^4\right\}^{1/2}$$

and

$$\left|E\left(\frac{1}{\sqrt{b}} \sum_{t=s}^{s+b-1} (X_t - EX_t)\right)^3\right|$$

$$\leq \left\{E\left(\frac{1}{\sqrt{b}} \sum_{t=s}^{s+b-1} (X_t - EX_t)\right)^4\right\}^{1/2} \left\{E\left(\frac{1}{\sqrt{b}} \sum_{t=s}^{s+b-1} (X_t - EX_t)\right)^2\right\}^{1/2}$$

$$\leq \left\{E\left(\frac{1}{\sqrt{b}} \sum_{t=s}^{s+b-1} (X_t - EX_t)\right)^4\right\}^{3/4},$$

by the Cauchy–Schwarz inequality. □

**Lemma A.6** (Leśkow and Synowiecki [21]). *Let $(X_1, \ldots, X_n)$ be a sample from the time series $\{X_t : t \in \mathbb{Z}\}$ that is real-valued and APC. Assume that the autocovariance function is uniformly summable, that is, there exists a summable sequence $\{c_\tau\}_{\tau=0}^{\infty}$ of real numbers such that $|\mathrm{Cov}(X_t, X_{t+\tau})| \leq c_\tau$. There then exists a number $\sigma^2$ such that, for any sequence $b = b(n) \leq n$ tending to infinity,*

$$\sup_{s=1,\ldots,n-b+1} \left| \mathrm{Var}\left( \frac{1}{\sqrt{b}} \sum_{t=s}^{s+b-1} X_t \right) - \sigma^2 \right| \longrightarrow 0 \quad \text{for } n \to \infty.$$

**Proof.** We refer the reader to Leśkow and Synowiecki [21]. □

**Proof of Theorem 3.1.** The following proof develops the techniques presented in Radulović [25] and Giné [11]. We can assume, without loss of generality, that $\mu = 0$; note that this does not mean that $EX_t \equiv 0$. For each $n \in \mathbb{N}$ and $t = 1, \ldots, n-b+1$ (recall that $b = b(n)$), we let $Z_{t,b} = X_t + \cdots + X_{t+b-1}$. The random variables $Z_{j,b}^*$ are conditionally independent (given the sample) with common distribution

$$P^*(Z_{j,b}^* = Z_{t,b}) = \frac{1}{n-b+1} \quad \text{for } t = 1, \ldots, n-b+1.$$

By Corollary 2.4.8 on page 63 of Araujo and Giné [1], the conclusion of the theorem is implied by the fact that for any $\delta > 0$,

$$\sum_{j=1}^{k} P^*\left( \frac{1}{\sqrt{n}} |Z_{j,b}^*| > \delta \right) \xrightarrow{P} 0, \tag{9}$$

$$\sum_{j=1}^{k} E^*\left( \frac{1}{\sqrt{n}} Z_{j,b}^* \mathbf{1}_{|Z_{j,b}^*| \leq \sqrt{n}\delta} \right) - \sum_{j=1}^{k} \frac{1}{\sqrt{n}} E^*(Z_{j,b}) \xrightarrow{P} 0 \tag{10}$$

and

$$\sum_{j=1}^{k} \mathrm{Var}^*\left( \frac{1}{\sqrt{n}} Z_{j,b}^* \mathbf{1}_{|Z_{j,b}^*| \leq \sqrt{n}\delta} \right) \xrightarrow{P} \sigma^2, \tag{11}$$

where $\mathbf{1}_S$ denotes the indicator function of the statement $S$.

In order to prove (9), observe that

$$\sum_{j=1}^{k} P^*\left( \frac{1}{\sqrt{n}} |Z_{j,b}^*| > \delta \right) = \frac{1}{n-b+1} \sum_{t=1}^{n-b+1} k \mathbf{1}_{|Z_{t,b}| > \sqrt{n}\delta}.$$

Let
$$U_{n,t} = k\mathbf{1}_{|Z_{t,b}|>\sqrt{n}\delta}$$

and consider the array $\{U_{n,t}: t = 1, \ldots, n - b + 1\}$. By Lemmas A.2 and A.3 with $f(x) = k\mathbf{1}_{|x|>C}$, we obtain that this array is row-wise strictly periodic. In order to show that $\{U_{n,t}\}$ satisfies the assumptions of Lemma A.4, define the triangular array $\{V_{n,t}: t = 1, \ldots, n - b + 1\}$ as

$$V_{n,t} = \frac{1}{b}Z_{t,b}^2.$$

Observe that for any fixed $t = 1, \ldots, T$,

$$\frac{1}{\sqrt{b}}Z_{t,b} = \frac{Z_{1,b+t-1}}{\sqrt{b+t-1}}\frac{\sqrt{b+t-1}}{\sqrt{b}} - \frac{Z_{1,t-1}}{\sqrt{b}} \xrightarrow{d} \mathcal{N}(\mu, \sigma^2)$$

and

$$E\left(\frac{1}{b}Z_{t,b}^2\right) = \operatorname{Var}\left(\frac{1}{\sqrt{b}}Z_{t,b}\right) + E^2\left(\frac{1}{\sqrt{b}}Z_{t,b}\right).$$

Since $\mu = 0$, the convergence $E^2((1/\sqrt{b})Z_{t,b}) \to 0$ is implied by the inequality (1), whereas the convergence $\operatorname{Var}((1/\sqrt{b})Z_{t,b}) \to \sigma^2$ results from Lemma A.6. Therefore, the sequences

$$\{V_{n,1}\}_{n=n_1}^{\infty}, \ldots, \{V_{n,T}\}_{n=n_T}^{\infty}$$

are uniformly integrable. Moreover, since $Z_{t,b}^2 > n\delta^2$ for $U_{n,t} \neq 0$, we obtain the estimation

$$0 \leq |U_{n,t}| = k\mathbf{1}_{|Z_{t,b}|>\sqrt{n}\delta} \leq \frac{Z_{t,b}^2}{\delta^2 b}\mathbf{1}_{|Z_{t,b}|>\sqrt{n}\delta} \leq \frac{1}{\delta^2}|V_{n,t}|,$$

which implies that the sequences

$$\{U_{n,1}\}_{n=n_1}^{\infty}, \ldots, \{U_{n,T}\}_{n=n_T}^{\infty}$$

are also uniformly integrable. Next, due to uniform integrability of $\{V_{n,t}\}$,

$$0 \leq E|U_{n,t}| = E(k\mathbf{1}_{|Z_{t,b}|>\sqrt{n}\delta})$$
$$\leq E\left(\frac{Z_{t,b}^2}{\delta^2 b}\mathbf{1}_{|Z_{t,b}|>\sqrt{n}\delta}\right) = \frac{1}{\delta^2}E|V_{n,t}|\mathbf{1}_{|V_{n,t}|>k\delta^2} \to 0,$$

so $EU_{n,t} \to 0$ for $n \to \infty$ and fixed $t = 1, \ldots, T$. It is easy to see that the triangular array $\{U_{n,t}\}$ is row-wise $\alpha$-mixing with $\alpha_{U_n}(\tau) = \alpha_X(\max\{\tau - b(n) + 1, 0\})$, where $\alpha_X$ is the mixing function of the underlying series $\{X_t\}$. By the inequality

$$\operatorname{Cov}(U_{n,t}\mathbf{1}_{|U_{n,t}|<A}, U_{n,t+\tau}\mathbf{1}_{|U_{n,t+\tau}|<A}) \leq 4A^2\alpha_{U_n}(\tau)$$

(Lemma A.0.2 from Politis *et al.* [24]) and the estimation

$$\frac{1}{n-b+1}\sum_{\tau=0}^{n-b}\alpha_{U_n}(\tau) \leq \frac{b}{n-b+1} + \frac{1}{n-b+1}\sum_{\tau=0}^{n-b}\alpha_X(\tau) \to 0 \qquad \text{for } n\to\infty,$$

we obtain that condition (iii) of Lemma A.4 is satisfied with $a_{n,\tau} = 4\alpha_{U_n}(\tau)$.

In order to prove (10), observe that

$$\sum_{j=1}^{k} E^*\left(\frac{1}{\sqrt{n}}Z_{j,b}^* \mathbf{1}_{|Z_{j,b}^*|\leq\sqrt{n}\delta}\right) - \sum_{j=1}^{k} \frac{1}{\sqrt{n}}E^*(Z_{j,b})$$

$$= \frac{1}{n-b+1}\sum_{t=1}^{n-b+1} k\frac{1}{\sqrt{n}}Z_{t,b}\mathbf{1}_{|Z_{t,b}|>\sqrt{n}\delta},$$

denote

$$U'_{n,t} = k\frac{1}{\sqrt{n}}Z_{t,b}\mathbf{1}_{|Z_{t,b}|>\sqrt{n}\delta}$$

and consider the array $\{U'_{n,t} : t=1,\ldots,n-b+1\}$. By Lemmas A.2 and A.3 with $f(x) = (k/\sqrt{n})x\mathbf{1}_{|x|>C}$, we have that this array is row-wise strictly periodic. Moreover, since $Z_{t,b}^2 > n\delta^2$ for $U'_{t,b} \neq 0$,

$$0 \leq |U'_{n,t}| = k\frac{1}{\sqrt{n}}|Z_{t,b}|\mathbf{1}_{|Z_{t,b}|>\sqrt{n}\delta} \leq k\frac{1}{n\delta}Z_{t,b}^2 = \frac{1}{\delta}\frac{Z_{t,b}^2}{b} = \frac{1}{\delta}V_{n,t},$$

so we have that the sequences

$$\{U'_{n,1}\}_{n=n_1}^{\infty}, \ldots, \{U'_{n,T}\}_{n=n_T}^{\infty}$$

are uniformly integrable. Since

$$0 \leq E|U'_{n,t}| = E\left(k\frac{1}{\sqrt{n}}|Z_{t,b}|\mathbf{1}_{|Z_{t,b}|>\sqrt{n}\delta}\right) \leq \frac{1}{\delta}E(|V_{n,t}|\mathbf{1}_{|V_{n,t}|>k\delta^2}) \to 0,$$

$EU'_{n,t} \to 0$ for $n\to\infty$ and fixed $t=1,\ldots,T$. By the same considerations as at the end of the proof of (9), we can apply Lemma A.4 to the array $\{U'_{n,t}\}$ to obtain that (10) is satisfied.

In order to prove (11), observe that

$$\sum_{j=1}^{k}\text{Var}^*\left(\frac{1}{\sqrt{n}}Z_{j,b}^*\mathbf{1}_{|Z_{j,b}^*|\leq\sqrt{n}\delta}\right)$$

$$= \frac{1}{n-b+1}\sum_{t=1}^{n-b+1}\frac{1}{b}Z_{t,b}^2\mathbf{1}_{|Z_{t,b}|\leq\sqrt{n}\delta} - \left(\frac{1}{n-b+1}\sum_{t=1}^{n-b+1}\frac{1}{\sqrt{b}}Z_{t,b}\mathbf{1}_{|Z_{t,b}|\leq\sqrt{n}\delta}\right)^2. \tag{12}$$

We will treat the terms separately. For the first, we have

$$\frac{1}{n-b+1}\sum_{t=1}^{n-b+1}\frac{1}{b}Z_{t,b}^2\mathbf{1}_{|Z_{t,b}|\leq\sqrt{n}\delta}$$

$$=\frac{1}{n-b+1}\sum_{t=1}^{n-b+1}\frac{1}{b}Z_{t,b}^2-\frac{1}{n-b+1}\sum_{t=1}^{n-b+1}\frac{1}{b}Z_{t,b}^2\mathbf{1}_{|Z_{t,b}|>\sqrt{n}\delta}$$

$$=\frac{1}{n-b+1}\sum_{t=1}^{n-b+1}V_{n,t}-\frac{1}{n-b+1}\sum_{t=1}^{n-b+1}V'_{n,t},$$

where

$$\{V'_{n,t}:t=1,\ldots,n-b+1\}=\left\{\frac{1}{b}Z_{t,b}^2\mathbf{1}_{|Z_{t,b}|>\sqrt{n}\delta}:t=1,\ldots,n-b+1\right\}.$$

We have shown uniform integrability of the sequences

$$\{V_{n,1}\}_{n=n_1}^{\infty},\ldots,\{V_{n,T}\}_{n=n_T}^{\infty}$$

and that for any fixed $t=1,\ldots,T$, $EV_{n,t}\to\sigma^2$. Due to this, and the fact that

$$V'_{n,t}=V_{n,t}\mathbf{1}_{|V_{n,t}|>k\delta^2},$$

we also have the convergence $EV'_{n,t}\to 0$ for fixed $t$ and $n\to\infty$. By Lemma A.4 applied to triangular arrays $\{V_{n,t}\}$ and $\{V'_{n,t}\}$, we obtain that the first term of the right-hand side of (12) tends in probability to $\sigma^2$. For the second term, we have

$$\frac{1}{n-b+1}\sum_{t=1}^{n-b+1}\frac{1}{\sqrt{b}}Z_{t,b}\mathbf{1}_{|Z_{t,b}|\leq\sqrt{n}\delta}$$

$$=\frac{1}{n-b+1}\sum_{t=1}^{n-b+1}\frac{1}{\sqrt{b}}Z_{t,b}-\frac{1}{n-b+1}\sum_{j=1}^{n-b+1}\frac{1}{\sqrt{b}}Z_{t,b}\mathbf{1}_{|Z_{t,b}|>\sqrt{n}\delta}.$$

Consider the arrays

$$\{T_{n,t}:t=1,\ldots,n-b+1\}=\left\{\frac{1}{\sqrt{b}}Z_{t,b}:t=1,\ldots,n-b+1\right\}$$

and

$$\{T'_{n,t}:t=1,\ldots,n-b+1\}=\left\{\frac{1}{\sqrt{b}}Z_{t,b}\mathbf{1}_{|Z_{t,b}|>\sqrt{n}\delta}:t=1,\ldots,n-b+1\right\}.$$

Since it is assumed that $\mu=0$ and that inequality (1) holds, we have

$$E(T_{n,t})=\frac{1}{\sqrt{b}}(X_t+\cdots+X_{t+b-1})=O\left(\frac{1}{\sqrt{b}}\right)\to 0\quad\text{for }n\to\infty.$$

By the Cauchy–Schwarz inequality, the sequences

$$\{T_{n,1}\}_{n=n_1}^{\infty}, \ldots, \{T_{n,T}\}_{n=n_T}^{\infty},$$
$$\{T'_{n,1}\}_{n=n_1}^{\infty}, \ldots, \{T'_{n,T}\}_{n=n_T}^{\infty}$$

are uniformly integrable. Applying Lemma A.4, we obtain that

$$\frac{1}{n-b+1} \sum_{t=1}^{n-b+1} \frac{1}{\sqrt{b}} Z_{t,b} \mathbf{1}_{|Z_{t,b}| \leq \sqrt{n}\delta} \xrightarrow{P} 0.$$

This completes the proof of (11). □

**Proof of Corollary 3.1.** We can write that

$$\frac{1}{n} \sum_{t=1}^{n} X_t = \frac{a_n}{\lfloor n/T \rfloor} \left\{ \frac{1}{T}\left(X_1 + \cdots + X_T\right) + \cdots \right.$$
$$\left. + \frac{1}{T}(X_{(\lfloor n/T \rfloor - 1)T + 1} + \cdots + X_{\lfloor n/T \rfloor T}) \right\} + o_P(1),$$

where

$$a_n = \frac{\lfloor n/T \rfloor T}{n} \to 1 \qquad \text{for } n \to \infty.$$

It is easy to see that the time series $Y_t = (X_{tT+1} + \cdots + X_{(t+1)T})/T$ is strictly stationary. Therefore, we may apply conditions for the CLT for the sample mean of strictly stationary sequences (Ibragimov and Linnik [16]) to see that under assumptions (i) and (ii), asymptotic normality holds. Now, note that, due to periodicity, we have

$$\sup_{t \in \mathbb{Z}} E|X_t|^{2+\delta} = M_{2+\delta} < \infty$$

and, by Lemma A.0.1 from Politis *et al.* [24],

$$\sum_{\tau=0}^{\infty} |\operatorname{Cov}(X_t, X_{t+\tau})| \leq 8 M_{2+\delta}^{2/(2+\delta)} \sum_{\tau=0}^{\infty} \alpha_X^{\delta/(2+\delta)}(\tau) < \infty.$$

By Theorem 3.1, the proof is completed. □

**Proof of Theorem 3.2.** Using the same steps as in the proof of Theorem 3.1, we must show the following laws of large numbers:

$$\frac{1}{n-b+1} \sum_{t=1}^{n-b+1} U_{n,t} \xrightarrow{P} 0;$$

$$\frac{1}{n-b+1}\sum_{t=1}^{n-b+1} U'_{n,t} \xrightarrow{P} 0;$$

$$\frac{1}{n-b+1}\sum_{t=1}^{n-b+1} V_{n,t} \xrightarrow{P} \sigma^2 + \mu^2;$$

$$\frac{1}{n-b+1}\sum_{t=1}^{n-b+1} V'_{n,t} \xrightarrow{P} 0;$$

$$\frac{1}{n-b+1}\sum_{t=1}^{n-b+1} T_{n,t} \xrightarrow{P} \mu;$$

$$\frac{1}{n-b+1}\sum_{t=1}^{n-b+1} T'_{n,t} \xrightarrow{P} 0,$$

where all of the arrays were defined in the proof of Theorem 3.1. Take $\{Y_{n,t}: t=1,\ldots,n-b+1\}$ to be any of these arrays. Due to the mixing condition of the underlying series $\{X_t\}$, we have the estimation

$$\mathrm{Cov}(Y_{n,t}\mathbf{1}_{|Y_{n,t}|<A}, Y_{n,t+\tau}\mathbf{1}_{|Y_{n,t+\tau}|<A}) \leq A^2 a_{n,\tau}$$

for $t=1,2,\ldots,\tau=0,1,\ldots,n-b+1-t$ and

$$\frac{1}{n-b+1}\sum_{\tau=0}^{n-b} a_{n,\tau} \to 0 \qquad \text{for } n\to\infty.$$

We must show that

$$\frac{1}{d_n}\sum_{t=1}^{d_n} Y_{n,t} \xrightarrow{P} L,$$

where $L$ is $0, \mu$ or $\sigma^2 + \mu^2$, depending on the array, and, in the sequel of the proof, without loss of generality, we assume that $\mu = 0$. Letting $\widetilde{Y}_{n,t} = Y_{n,t}\mathbf{1}_{|Y_{n,t}|<A_n}$ and $\overline{Y}_{n,t} = Y_{n,t}\mathbf{1}_{|Y_{n,t}|\geq A_n}$, the above convergence follows from

$$\frac{1}{d_n}\sum_{t=1}^{d_n}(\widetilde{Y}_{n,t} - E\widetilde{Y}_{n,t}) \xrightarrow{P} 0, \tag{13}$$

$$\frac{1}{d_n}\sum_{t=1}^{d_n}(\overline{Y}_{n,t} - E\overline{Y}_{n,t}) \xrightarrow{P} 0 \tag{14}$$

and

$$\frac{1}{d_n}\sum_{t=1}^{d_n} EY_{n,t} \to L. \tag{15}$$

Taking $A_n$ as in the proof of Lemma A.4, condition (14) is satisfied due to the same considerations as in the proof of (7). Therefore, it suffices to show that conditions (14) and (15) hold.

We have that, for any $\epsilon > 0$,

$$P\left(\left|\frac{1}{d_n}\sum_{t=1}^{d_n}(\overline{Y}_{n,t} - E\overline{Y}_{n,t})\right| > \epsilon\right) \leq \frac{1}{d_n\epsilon}E\left|\sum_{t=1}^{d_n}(\overline{Y}_{n,t} - E\overline{Y}_{n,t})\right|$$

$$\leq \frac{1}{d_n\epsilon}\sum_{t=1}^{d_n}E|\overline{Y}_{n,t} - E\overline{Y}_{n,t}|$$

$$\leq \frac{2}{\epsilon}\left(\frac{1}{d_n}\sum_{t=1}^{d_n}E|\overline{Y}_{n,t}|\right).$$

In the following, we will carefully use the Cauchy–Schwarz and Chebyshev inequalities. For the array $\{U_{n,t}\}$, we have

$$\frac{1}{n-b+1}\sum_{t=1}^{n-b+1}E(|U_{n,t}|\mathbf{1}_{|U_{n,t}|\geq A_n}) \leq \frac{1}{n-b+1}\sum_{t=1}^{n-b+1}E(k\mathbf{1}_{|Z_{t,b}|>\sqrt{n}\delta})$$

$$= \frac{1}{n-b+1}\sum_{t=1}^{n-b+1}kP(|Z_{t,b}| > \sqrt{n}\delta)$$

$$\leq \frac{1}{n-b+1}\sum_{t=1}^{n-b+1}k\frac{E|Z_{t,b}|^4}{n^2\delta^4}$$

$$= \frac{1}{k\delta^4}\left\{\frac{1}{n-b+1}\sum_{t=1}^{n-b+1}E\left(\frac{1}{b^2}Z_{t,b}^4\right)\right\}.$$

For the array $\{U'_{n,t}\}$, we have

$$\frac{1}{n-b+1}\sum_{t=1}^{n-b+1}E(|U'_{n,t}|\mathbf{1}_{|U'_{n,t}|\geq A_n}) \leq \frac{1}{n-b+1}\sum_{t=1}^{n-b+1}E\left(\frac{k}{\sqrt{n}}|Z_{t,b}|\mathbf{1}_{|Z_{t,b}|>\sqrt{n}\delta}\right)$$

$$\leq \frac{1}{n-b+1}\sum_{t=1}^{n-b+1}\sqrt{k}E^{1/2}\left(\frac{1}{b}Z_{t,b}^2\right)P^{1/2}(|Z_{t,b}| > \sqrt{n}\delta)$$

$$\leq \frac{1}{\sqrt{k}\delta^2}\left\{\frac{1}{n-b+1}\sum_{t=1}^{n-b+1}E^{3/4}\left(\frac{1}{b^2}Z_{t,b}^4\right)\right\}.$$

For the array $\{V_{n,t}\}$, we have

$$\frac{1}{n-b+1}\sum_{t=1}^{n-b+1}E(|V_{n,t}|\mathbf{1}_{|V_{n,t}|\geq A_n}) = \frac{1}{n-b+1}\sum_{t=1}^{n-b+1}E\left(\frac{1}{b}Z_{t,b}^2\mathbf{1}_{(1/b)Z_{t,b}^2\geq A_n}\right)$$

$$\leq \frac{1}{n-b+1}\sum_{t=1}^{n-b+1}E^{1/2}\left(\frac{1}{b^2}Z_{t,b}^4\right)P^{1/2}\left(\frac{1}{b}Z_{t,b}^2\geq A_n\right)$$

$$\leq \frac{1}{n-b+1}\sum_{t=1}^{n-b+1}E^{1/2}\left(\frac{1}{b^2}Z_{t,b}^4\right)\frac{E^{1/2}((1/b^2)Z_{t,b}^4)}{A_n}$$

$$\leq \frac{1}{A_n}\left\{\frac{1}{n-b+1}\sum_{t=1}^{n-b+1}E\left(\frac{1}{b^2}Z_{t,b}^4\right)\right\}.$$

For the array $\{V'_{n,t}\}$, we have

$$\frac{1}{n-b+1}\sum_{t=1}^{n-b+1}E(|V'_{n,t}|\mathbf{1}_{|V'_{n,t}|\geq A_n}) \leq \frac{1}{n-b+1}\sum_{t=1}^{n-b+1}E\left(\frac{1}{b}Z_{t,b}^2\mathbf{1}_{|Z_{t,b}|>\sqrt{n}\delta}\right)$$

$$\leq \frac{1}{n-b+1}\sum_{t=1}^{n-b+1}E^{1/2}\left(\frac{1}{b^2}Z_{t,b}^4\right)P^{1/2}(|Z_{t,b}|>\sqrt{n}\delta)$$

$$\leq \frac{1}{k\delta^2}\left\{\frac{1}{n-b+1}\sum_{t=1}^{n-b+1}E\left(\frac{1}{b^2}Z_{t,b}^4\right)\right\}.$$

For the array $\{T_{n,t}\}$, we have

$$\frac{1}{n-b+1}\sum_{t=1}^{n-b+1}E(|T_{n,t}|\mathbf{1}_{|T_{n,t}|\geq A_n}) = \frac{1}{n-b+1}\sum_{t=1}^{n-b+1}E\left(\frac{1}{\sqrt{b}}|Z_{t,b}|\mathbf{1}_{(1/\sqrt{b})|Z_{t,b}|\geq A_n}\right)$$

$$\leq \frac{1}{n-b+1}\sum_{t=1}^{n-b+1}E^{1/2}\left(\frac{1}{b}Z_{t,b}^2\right)P^{1/2}\left(\frac{1}{\sqrt{b}}|Z_{t,b}|\geq A_n\right)$$

$$\leq \frac{1}{A_n^2}\left\{\frac{1}{n-b+1}\sum_{t=1}^{n-b+1}E^{3/4}\left(\frac{1}{b^2}Z_{t,b}^4\right)\right\}.$$

For the array $\{T'_{n,t}\}$, we have

$$\frac{1}{n-b+1}\sum_{t=1}^{n-b+1} E(|T'_{n,t}|\mathbf{1}_{|T'_{n,t}|\geq A_n}) \leq \frac{1}{n-b+1}\sum_{t=1}^{n-b+1} E\left(\frac{1}{\sqrt{b}}|Z_{t,b}|\mathbf{1}_{|Z_{t,b}|>\sqrt{n}\delta}\right)$$

$$\leq \frac{1}{n-b+1}\sum_{t=1}^{n-b+1} E^{1/2}\left(\frac{1}{b}Z_{t,b}^2\right)P^{1/2}(|Z_{t,b}|>\sqrt{n}\delta)$$

$$\leq \frac{1}{k\delta^2}\left\{\frac{1}{n-b+1}\sum_{t=1}^{n-b+1} E^{3/4}\left(\frac{1}{b^2}Z_{t,b}^4\right)\right\}.$$

Applying assumption (iii), Lemma A.5 and the fact that both $k\to\infty$ and $A_n\to\infty$, we obtain that condition (14) is satisfied for all of the arrays. As for condition (15), from the above, we immediately obtain that for $n\to\infty$,

$$\frac{1}{n-b+1}\sum_{t=1}^{n-b+1} EU_{n,t} \to 0,$$

$$\frac{1}{n-b+1}\sum_{t=1}^{n-b+1} EU'_{n,t} \to 0,$$

$$\frac{1}{n-b+1}\sum_{t=1}^{n-b+1} EV'_{n,t} \to 0,$$

$$\frac{1}{n-b+1}\sum_{t=1}^{n-b+1} ET'_{n,t} \to 0.$$

Moreover,

$$\frac{1}{n-b+1}\sum_{t=1}^{n-b+1} ET_{n,t} = \frac{1}{n-b+1}\sum_{t=1}^{n-b+1} \frac{1}{\sqrt{b}}EZ_{t,b} \to 0,$$

which follows from inequality (1). As for the remaining array $\{V_{n,t}\}$, we have

$$\frac{1}{n-b+1}\sum_{t=1}^{n-b+1} EV_{n,t} = \frac{1}{n-b+1}\sum_{t=1}^{n-b+1} \frac{1}{b}EZ_{t,b}^2$$

$$= \frac{1}{n-b+1}\sum_{t=1}^{n-b+1}\left\{\mathrm{Var}\left(\frac{1}{\sqrt{b}}(X_t+\cdots+X_{t+b-1})\right)\right.$$

$$\left.+\frac{1}{b}(EX_t+\cdots+EX_{t+b-1})^2\right\} \to \sigma^2,$$

which is implied by Lemma A.6 and inequality (1). This remark completes the proof of Theorem 3.2. $\square$

**Proof of Corollaries 3.2 and 3.3.** We use Theorem 3.2. By the same reasoning as in the proof of Corollary 3.1, we obtain that the autocovariance function is uniformly summable. Assumption (iii) of Theorem 3.2 is implied by Theorem 5 from Kim [17] for Corollary 3.2 and by the reasoning of Annexe C from Rio [26] for Corollary 3.3. $\square$

**Proof of Theorem 4.1.** This is straightforwardly implied by the Cramér–Wold device, Corollary 3.2 and Remark 4.1. We simply note than any linear combination of the real and imaginary parts of the series $\{W_t(\lambda,\tau)\}$ also satisfies the assumptions of Corollary 3.2. $\square$

# Acknowledgements

The author is grateful to Prof. J. Leśkow for reading an early version of this manuscript. This work was partially supported by AGH local grant 10.420.03.

# References

[1] Araujo, A. and Giné, E. (1980). *The Central Limit Theorem for Real and Banach Valued Random Variables*. New York: Wiley. MR0576407
[2] Arcones, M. and Giné, E. (1989). The bootstrap of the mean with arbitrary bootstrap sample size. *Ann. Inst. H. Poincaré Probab. Statist.* **25** 457–481. MR1045246
[3] Bloomfield, P., Hurd, H. and Lund, R. (1994). Periodic correlation in stratospheric ozone time series. *J. Time Ser. Anal.* **15** 127–150. MR1263886
[4] Cambanis, S., Houdré, C., Hurd, H. and Leśkow, J. (1994). Laws of large numbers for periodically and almost periodically correlated processes. *Stochastic Process. Appl.* **53** 37–54. MR1290706
[5] Chan, V., Lahiri, S.N. and Meer, W. (2004). Block bootstrap estimation of the distribution of cumulative outdoor degradation. *Technometrics* **46** 215–224. MR2060017
[6] Corduneanu, C. (1989). *Almost Periodic Functions*. New York: Chelsea.
[7] Dehay, D. and Leśkow, J. (1996). Functional limit theory for the spectral covariance estimator. *J. Appl. Probab.* **33** 1077–1092. MR1416228
[8] Dehay, D. and Leśkow, J. (1996). Testing stationarity for stock market data. *Econom. Lett.* **50** 205–212.
[9] Fitzenberger, B. (1997). The moving blocks bootstrap and robust inference for linear last squares and quantile regressions. *J. Econometrics* **82** 235–287. MR1613422
[10] Gardner, W., Napolitano, A. and Paura, L. (2006). Cyclostationarity: Half a century of research. *Signal Processing* **86** 639–697.
[11] Giné, E. (1997). Lectures on some aspects of the bootstrap. In *Lectures on Probability Theory and Statistics* (*Saint-Flour, 1996*). *Lecture Notes in Math.* **1665** 37–151. Berlin: Springer. MR1490044
[12] Gladyshev, E.G. (1961). Periodically correlated random sequences. *Sov. Math.* **2** 385–388.
[13] Gonçalves, S. and White, H. (2002). The bootstrap of the mean for dependent heterogeneous arrays. *Econometric Theory* **18** 1367–1384. MR1945417

[14] Hurd, H. (1991). Correlation theory of almost periodically correlated processes. *J. Multivariate Anal.* **30** 24–45. MR1097303
[15] Hurd, H. and Leśkow, J. (1992). Strongly consistent and asymptotically normal estimation of the covariance for almost periodically correlated processes. *Statist. Decisions* **10** 201–225. MR1183203
[16] Ibragimov, I.A. and Linnik, Y.V. (1971). *Independent and Stationary Sequences of Random Variables*. Groningen: Wolters-Noordhoff. MR0322926
[17] Kim, T.Y. (1994). Moment bounds for nonstationary dependent sequences. *J. Appl. Probab.* **31** 731–742. MR1285511
[18] Künsch, H. (1989). The jackknife and the bootstrap for general stationary observations. *Ann. Statist.* **17** 1217–1241. MR1015147
[19] Lahiri, S.N. (1992). Edgeworth correction by moving block bootstrap for stationary and nonstationary data. In *Exploring the Limits of Bootstrap* (R. LePage and L. Billard, eds.) 183–214. New York: Wiley. MR1197785
[20] Lahiri, S.N. (1999). Theoretical comparison of block bootstrap methods. *Ann. Statist.* **27** 386–404. MR1701117
[21] Leśkow, J. and Synowiecki, R. (2006). Asymptotic distribution and subsampling for the estimator of autocovariance coefficient for APC time series. Preprint No. 1/2006, Faculty of Applied Mathematics, AGH University of Science and Technology.
[22] Liu, R.Y. and Singh, K. (1992). Moving blocks jacknife and bootstrap capture weak dependence. In *Exploring the Limits of Bootstrap* (R. LePage and L. Billard, eds.) 225–248. New York: Wiley. MR1197787
[23] Politis, D. (2001). Resampling time series with seasonal components. In *Proceedings of the 33rd Symposium on the Interface of Computing Science and Statistics, Orange County, California, June 13–17*, 619–621.
[24] Politis, D., Romano, J. and Wolf, M. (1999). *Subsampling*. New York: Springer. MR1707286
[25] Radulović, D. (1996). The bootstrap of the mean for strong mixing sequences under minimal conditions. *Statist. Probab. Lett.* **28** 65–72. MR1394420
[26] Rio, E. (2000). *Theorie asymptotique des processus aleatoires faiblement dependants*. Berlin: Springer. MR2117923
[27] Synowiecki, R. (2007). Some results on the subsampling for $\varphi$-mixing periodically strictly stationary time series. *Probab. Math. Statist.* **27**. To appear.
[28] Yeung, G.K. and Gardner, W. (1996). Search-efficient methods of detection of cyclostationarity signals. *IEEE Trans. Signal Process.* **44** 1214–1223.